\documentclass[12pt,scs]{article}

\newtheorem{thm}{Theorem}[section]

\newtheorem{lem}[thm]{Lemma}

\begin{document}

\title{{\bf A Smooth and Compactly Supported Radial Function}}
\author{{\bf Lin-Tian Luh\thanks{This work was supported by NSC93-2115-M126-004.}}\\
Department of Mathematics, Providence University \\Shalu, Taichung\\
email:ltluh@pu.edu.tw\\phone:(04)26328001 ext. 15126\\ fax:(04)26324653}

\maketitle

\bigskip
{\bf Abstract.} In the field of radial basis functions mathematicians have been endeavouring to find infinitely differentiable and compactly supported radial functions. This kind of functions are extremely important for some reasons. First, its computational properties will be very good since it's compactly supported. Second, its error bound will converge very fast since it's infinitely differentiable. However there is hitherto no such functions which can be expressed in a simple form. This is a famous question. The purpose of this paper is to answer this question.\\
\\
{\bf keywords}: radial basis function, compactly supported function, positive definite function.

\pagenumbering{arabic}
\setcounter{section}{0}
\section{Introduction}
In this paper we raise a new function
\begin{equation}
  \Phi(x):=\left\{ \begin{array}{ll}
                     e^{\frac{-\alpha}{1-\| x\| ^{2}}} & \mbox{if $\| x\| <1$}\\
                     0                                 & \mbox{if $\| x\| \geq 1$}
                   \end{array}
           \right.  
\end{equation} 
where $\alpha>0$ is a constant and $\| x\| $ denotes the Euclidean norm of $x\in R^{n}$. The infinite differentiability of $\Phi$ can be shown by routine investigation, and we omit it. What's difficult is to show its positive definiteness. In order to make it useful in scattered data interpolation, this problem has to be overcome. This is the main theme of the next section.
\section{Positive Definiteness}
It's well known that positive definiteness is just conditional positive definiteness of order zero. Here we adopt the definition of [1] for conditional positive definiteness. Before showing that the map defined in (1) is positive definite, let's define
\begin{equation}
  \phi(t):=\left\{ \begin{array}{ll}
                     e^{\frac{-\alpha}{1-t^{2}}} & \mbox{if $0\leq t<1$}\\
                     0                           & \mbox{if $1\leq t.$}
                   \end{array}
           \right.  
\end{equation}
where $\alpha>0$ is a constant. In other words, $\phi(\| x\| )=\Phi(x)$ for $x\in R^{n}$. In order to make this paper easier to understand, let's analyze the derivatives of $\phi(\sqrt{t})$ of order up to eight.

It's easily seen that with $u=\frac{1}{1-t}$, we have $\frac{du}{dt}=u^{2}$ for $u>1$ and
\begin{itemize}
  \item $\frac{d}{dt}\phi(\sqrt{t})=e^{-\alpha u}(-\alpha)u^{2}$
  \item $\frac{d^{2}}{dt^{2}}\phi(\sqrt{t})=e^{-\alpha u}[\alpha^{2}u^{4}-2\alpha u^{3}]$
  \item $\frac{d^{3}}{dt^{3}}\phi(\sqrt{t})=e^{-\alpha u}[-\alpha^{3}u^{6}+6\alpha^{2}u^{5}-6\alpha u^{4}]$
  \item $\frac{d^{4}}{dt^{4}}\phi(\sqrt{t})=e^{-\alpha u}[\alpha^{4}u^{8}-12\alpha^{3}u^{7}+36\alpha^{2}u^{6}-24\alpha u^{5}]$
  \item $\frac{d^{5}}{dt^{5}}\phi(\sqrt{t})=e^{-\alpha u}[-\alpha^{5}u^{10}+20\alpha^{4}u^{9}-120\alpha^{3}u^{8}+240\alpha^{2}u^{7}-120\alpha u^{6}]$
  \item $\frac{d^{6}}{dt^{6}}\phi(\sqrt{t})=e^{-\alpha u}[\alpha^{6}u^{12}-30\alpha^{5}u^{11}+300\alpha^{4}u^{10}-1200\alpha^{3}u^{9}+1800\alpha^{2}u^{8}-720\alpha u^{7}]$
  \item $\frac{d^{7}}{dt^{7}}=e^{-\alpha u}[-\alpha^{7}u^{14}+42\alpha^{6}u^{13}-630\alpha^{5}u^{12}+4200\alpha^{4}u^{11}-12600\alpha^{3}u^{10}+15120\alpha^{2}u^{9}-5040\alpha u^{8}]$
  \item $\frac{d^{8}}{dt^{8}}\phi(\sqrt{t})=e^{-\alpha u}[\alpha^{8}u^{16}-56\alpha^{7}u^{15}+1176\alpha^{6}u^{14}-11760\alpha^{5}u^{13}+58800\alpha^{4}u^{12}-141120\alpha^{3}u^{11}+141120\alpha^{2}u^{10}-40320\alpha u^{9}].$
\end{itemize}
By induction, for $j\geq 1$,
$$\frac{d^{j}}{dt^{j}}\phi(\sqrt{t})=e^{-\alpha u}[(-1)^{j}\alpha^{j}u^{2j}+(-1)^{j-1}a_{j-1}\alpha^{j-1}u^{2j-1}+\cdots +(-1)a_{1}\alpha u^{j+1}]$$
where $a_{1},\cdots ,a_{j-1}$ are nonnegative integers. Obviously
$$(-1)^{j}\frac{d^{j}}{dt^{j}}\phi(\sqrt{t})\geq 0$$
for all $j\geq 1$ and $u\geq 1 $(i.e. $0\leq t<1$) as long as $\alpha$ is large enough. In fact we have the following more accurate lemma.
\begin{lem}
  Assume that $\alpha_{1}=0$ and $\alpha_{j+1}=\alpha_{j}+2j$ for $j=1,2,3,\cdots $. Then for any $j\geq 1$, 
$$(-1)^{j}\frac{d^{j}}{dt^{j}}\phi(\sqrt{t})\geq 0$$
for all $t>0$ if $\alpha \geq \alpha_{j}$.
\end{lem} 
Proof. Note that if $\alpha_{j}$ is defined as in the lemma, then $\alpha_{j}$ is just the exponent of $u$ in the first term of $\frac{d^{j-1}}{dt^{j-1}}\phi(\sqrt{t})$ plus the absolute value of the integer coefficient in the second term of $\frac{d^{j-1}}{dt^{j-1}}\phi(\sqrt{t})$. In turn the absolute value of the integer coefficient in the second term of $\frac{d^{j-1}}{dt^{j-1}}\phi(\sqrt{t})$ is the exponent of $u$ in the first term of $\frac{d^{j-2}}{dt^{j-2}}\phi(\sqrt{t})$ plus the absolute value of the integer coefficient in the second term of $\frac{d^{j-2}}{dt^{j-2}}\phi (\sqrt{t})$. Our lemma thus follows immediately by simple induction.    
\begin{thm}
  (a)Suppose $\alpha\geq 12$. Then $\Phi(x)$ in (1) is positive definite on $R^{d}$ for $1\leq d \leq 4$. (b)Suppose $n>4$ and $\alpha \geq \alpha_{j}$ for $j=\lceil \frac{n-4}{2}\rceil +4$, where $\lceil \lambda \rceil $ denotes the least integer not less than $\lambda$. Then $\Phi(x)$ in (1) is positive definite on $R^{d}$ for $1\leq d\leq n$.
\end{thm}
Proof. This is an immediate result of Theorem1 of [1], Theorem2.2 of [2] and Lemma2.1 of this paper. \\
\\
{\bf Remark}: In this theorem we treat $R^{d}$ for $1\leq d\leq 4$ independently to make it easy for use. In fact, if we require $\Phi$ be positive definite on $R^{d}$ for $1\leq d\leq 2$ only, then $\alpha \geq 6$ is sufficient. The first values of $\alpha_{j}$ are $0,2,6,12,20,30,42,56,72,90,\cdots$ ; for $j=1,2,3,4,5,6,7,8,9,10,\cdots$, respectively. It's easily seen from our deduction of $\frac{d^{j}}{dt^{j}}\phi(\sqrt{t})$ that $$(-1)^{j}\frac{d^{j}}{dt^{j}}\phi(\sqrt{t})\geq 0$$ for $j=1,2,\cdots ,8$ if $\alpha\geq \alpha_{j}$, and $t>0$. Thus we already checked that $\Phi(x)$ is positive definite on $R^{d}$ for $1\leq d\leq 12$ if $\alpha \geq 56$. For higher dimensions, the values of the necessary $\alpha_{j}$ can be calculated recursively. It's a bit more tedious, but very easy.\\


\begin{thebibliography}{99}
\bibitem{Dy}N. Dyn,
{\em Interpolation and Approximation by Radial and Related Functions,}
Approximation Theory VI, (C.K. Chui, L.L. Schumaker and J. Ward eds.), Academic Press,(1989),211-234.

\bibitem{Mi}C.A. Micchelli,
{\em Interpolation of Scattered Data:Distance Matrices and Conditionally Positive Definite Functions,}
Constr. Approx. 2(1986),11-22. 
\end{thebibliography}
\end{document}